\documentstyle[11pt,leqno]{article}
\newtheorem{theorem}{{\sc Theorem}}
\newcommand{\bt}{\begin{theorem}}
\newcommand{\et}{\end{theorem}}
\newcommand{\newsection}[1]{\setcounter{equation}{0} \setcounter{theorem}{0}
\section{#1}}

\newcommand{\NI}{\noindent}
\newcommand{\bea}{\begin{eqnarray}}
\newcommand{\eea}{\end{eqnarray}}

\def \spec#1 {\mathop{#1}}

\def \b #1 {\bf #1}
\newcommand {\nnb }{\nonumber}

\newcommand{\clf}{{\cal F}}

\newcommand{\ity}{\infty}
\newcommand{\raro}{\rightarrow}

\newcommand{\vsp}{\vskip 1em}

\newcommand{\be}{\begin{equation}}
\newcommand{\ee}{\end{equation}}
\newcommand{\ben}{\begin{eqnarray*}}
\newcommand{\een}{\end{eqnarray*}}

\begin{document}
\noindent{Research article}\\
\vsp
\noindent{B.L.S. Prakasa Rao{\footnote{{\bf Corresponding author:  B.L.S. Prakasa Rao}, CR RAO Advanced Institute of Mathematics, Statistics and Computer Science, Hyderabad, India, e-mail:blsprao@gmail.com}}}\\
\vsp
\noindent{\Large{\bf Maximum likelihood estimation for}}\\
\noindent{\Large{\bf sub-fractional Vasicek Model}}\\
\vsp
\NI{\bf Abstract:} We investigate the asymptotic properties of maximum likelihood estimators of the drift parameters  for fractional Vasicek model driven by a sub-fractional Brownian motion.\\
\vsp
\NI{\bf Keywords and phrases}: Sub-fractional Vasicek model ; sub-fractional Brownian motion; maximum likelihood estimation.
\vsp
\NI {\bf MSC 2020:} Primary 62M09, Secondary 60G22.
\vsp
\newsection{Introduction}

Statistical inference  for fractional diffusion processes satisfying stochastic
differential equations driven by a fractional Brownian motion (fBm)  has been studied
earlier and a comprehensive survey of various methods is given in Prakasa
Rao [21] and Mishura [18]. Processes driven by fBm with Hurst index $H\frac{1}{2}$ have been used for modeling purposes whenever there is long range dependence (cf. Prakasa Rao [21]). However in some applications such as turbulence phenomena in hydromechanics, it was found that  fbm is adequate for modeling small increments but it seems to be inadequate for large increments. For this reason, a sub-fractional Brownian motion may be an alternative to fBm for modeling (cf. Mishura and Zili [17]).  There has been a recent interest to study inference problems for stochastic processes driven by a sub-fractional Brownian motion. Bojdecki et al. [1] introduced a centered Gaussian process $\zeta^H= \{\zeta^H(t), t\geq 0\}$ called  {\it sub-fractional Brownian motion} (sub-fBm) with the covariance function
$$C_H(s,t)= s^{2H}+t^{2H}-\frac{1}{2}[(s+t)^{2H}+|s-t|^{2H}]$$
where $0<H<1.$ The increments of this process are not stationary and are more weakly correlated on non-overlapping intervals than those of a fBm. Tudor [34] introduced a Wiener integral with respect to a sub-fBm. Tudor [31,32,33,34] discussed some properties related to sub-fBm and its corresponding stochastic calculus. By using a fundamental martingale associated to sub-fBm, a Girsanov type theorem is obtained in Tudor [34]. Diedhiou et al. [4] investigated parametric estimation for a stochastic differential equation (SDE) driven by a sub-fBm. Mendy [16] studied parameter estimation for the sub-fractional Ornstein-Uhlenbeck process defined by the stochastic differential equation
$$dX_t=\theta X_tdt+d\zeta^H(t), t \geq 0$$ where $H>\frac{1}{2}.$ This is an analogue of the Ornstein-Uhlenbeck process, that is, a continuous time first order autoregressive process $X=\{X_t, t \geq 0\}$ which is the solution of a one-dimensional homogeneous linear stochastic differential equation driven  by a sub-fBm $\zeta^H= \{\zeta_t^H, t \geq 0\}$ with Hurst parameter $H.$ Mendy [16] proved that the least squares estimator estimator $\tilde \theta_T$ is strongly consistent as $T \raro \ity.$
Kuang and Xie [10] studied properties of maximum likelihood estimator for sub-fBm through approximation by a random walk. Kuang and Liu [9]discussed about the $L^2$-consistency  and strong consistency of the maximum likelihood estimators for the sub-fBm with drift based on discrete observations. Yan et al. [40] obtained the Ito's formula for sub-fractional Brownian motion with Hurst index $H>\frac{1}{2}.$  Shen and Yan [41] studied estimation for the drift of sub-fractional Brownian motion and constructed a class of biased estimators of James-Stein type which dominate the maximum likelihood estimator under the quadratic risk. El Machkouri et al. [6] investigated the asymptotic properties of the least squares estimator for non-ergodic Ornstein-Uhlenbeck process driven by Gaussian processes, in particular, sub-fractional Brownian motion. Es-sebaiy and Es-sebaiy [7] investigated the problem of estimation of drift parameters in a non-ergodic fractional vasicek model. In a recent paper, we have investigated optimal estimation of a signal perturbed by a sub-fractional Brownian motion in Prakasa Rao [23]. Some maximal and integral inequalities for a sub-fBm were derived in Prakasa Rao [22,26].  Parametric estimation for linear stochastic differential equations driven by a sub-fractional Brownian motion is studied in Prakasa Rao [24].  A Berry-Esseen type bound for the distribution of the maximum likelihood estimator for the drift parameter of a fractional Ornstein-Uhlenbeck type process driven by a sub-fractional Brownian motion is studied in Prakasa Rao [25]. Nonparametric estimation of trend for stochastic differential equations driven by sub-fractional Brownian motion is investigated in Prakasa Rao [27]. Mishura and Zili [17] gives a survey of stochastic analysis of mixed fractional Gaussian processes including a sub-fractional Brownian motion.

Vasicek [35] introduced an interest rate model
$$dX_t=(\alpha-\beta X_t)dt+\gamma dW_t, X_0=x_0\in R$$
where $\{W_t,t\geq 0\}$ is the standard Brownian motion for studying the market behaviour due to changes in the interest rates. This model is now known as the Vasicek model and it was generalized to a fractional Vasicek model to study processes with long range dependence which appears in financial mathematics and other areas such as telecommunication networks, turbulence and image processing. Properties of fractional Vasicek model for modeling are investigated  in Chronopoulu and Viens [2], Corlay et al. [3], Hao et al. [8], Song and Li [30] and Xiao et al. [37] among others. Maximum likelihood estimation for fractional Vasicek model is investigated in Lohvinenko and Ralchenko [12,13,14,15] and Xiao and Yu [38,39]. Maximum likelihood estimation for a fractional Vasicek model driven by a mixed fractional Brownian motion is studied in Prakasa Rao [28]. 

Our aim in this paper is to investigate the problem of maximum likelihood estimation of parameters in a  sub-fractional Vasicek model for processes driven by a sub-fractional Brownian motion.

\newsection{Sub-fractional Brownian motion}

Let $(\Omega, \clf, (\clf_t), P) $ be a stochastic basis satisfying the
usual conditions and the processes discussed in the following are
$(\clf_t)$-adapted. Further the natural filtration of a process is
understood as the $P$-completion of the filtration generated by this
process.

Let $\zeta^H= \{\zeta_t^H, t \geq 0 \}$ be a  normalized  {\it sub-fractional Brownian
motion} (sub-fBm)  with Hurst parameter $H \in (0,1)$, that is, a Gaussian process
with continuous sample paths such that $\zeta_0^H=0, E(\zeta_t^H)=0$ and
\be
E(\zeta_s^H \zeta_t^H)= t^{2H}+s^{2H}-\frac{1}{2}[(s+t)^{2H}+|s-t|^{2H}], t \geq 0, s \geq 0.
\ee

Bojdecki et al. [1] noted that the process
$$\frac{1}{\sqrt{2}}[W^H(t)+W^H(-t)], t \geq 0,$$
where $\{W^H(t), -\ity<t<\ity\}$ is a fBm, is a centered Gaussian process with the same covariance function as that of a  sub-fBm. This proves the existence of a sub-fBm. Let $D_H(s,t)$ denote the covariance function of a standard fractional Brownian motion with Hurst index $H.$  Note that
$$D_H(s,t)= \frac{1}{2}(|t|^{2H}+|s|^{2H}-|t-s|^{2H}).$$
Bojdecki et al. [1] proved the following result concerning properties of a sub-fBm.
\vsp
\NI{\bf Theorem 2.1:} {\it Let $\zeta^H= \{\zeta^H(t), t \geq 0\}$ be a sub-fBm defined on a filtered probability space $(\Omega, \clf, (\clf_t,t\geq 0),P).$ Then the following properties hold.

(i) The process $\zeta^H$ is self-similar, that is, for every $a>0,$
$$\{\zeta^H(at),t\geq 0\} \stackrel {\Delta}=  \{a^H\zeta^H(t), t\geq 0\}$$
in the sense that the processes, on both the sides of the equality sign, have the same finite dimensional distributions.

(ii) The process $\zeta^H$ is not Markov and it is not a semi-martingale.

(iii) For all $s,t \geq 0,$ the covariance function $C_H(s,t)$ of the process $\zeta^H$ is positive for all $s>0,t>0.$ Furthermore

$$C_H(s,t)>D_H(s,t) \;\;\mbox{if}\;\;H<\frac{1}{2}$$
and
$$C_H(s,t)< D_H(s,t) \;\;\mbox{if}\;\;H>\frac{1}{2}.$$

(iv) Let $\beta_H= 2-2^{2H-1}.$ For all $s\geq 0, t\geq 0,$

$$\beta_H(t-s)^{2H}\leq E[\zeta^H(t)-\zeta^H(s)]^2\leq (t-s)^{2H},\;\;\mbox{if}\;\;H>\frac{1}{2}$$
and
$$(t-s)^{2H}\leq E[\zeta^H(t)-\zeta^H(s)]^2\leq \beta_H (t-s)^{2H},\;\;\mbox{if}\;\;H<\frac{1}{2}$$
and the constants in the above inequalities are sharp.

(v) The process $\zeta^H $ has continuous sample paths almost surely and, for each $0<\epsilon<H$ and $T>0,$  there exists a random variable $K_{\epsilon,T}$ such that}
$$|\zeta^H(t)-\zeta^H(s)|\leq K_{\epsilon,T} |t-s|^{H-\epsilon}, 0\leq s,t \leq T.$$
\vsp
Let $f:[0,T]\raro R$ be a measurable  function and $\alpha  >0,$ and  $\sigma$ and $\eta$ be real. Define the Erdeyli-Kober-type fractional integral
\be
(I_{T,\sigma,\eta}^\alpha f)(s)=\frac{\sigma s^{\sigma \eta}}{\Gamma(\alpha)}\int_s^T\frac{t^{\sigma(1-\alpha-\eta)-1}f(t)}{(t^\sigma-s^\sigma)^{1-\alpha}}dt, s\in [0,T],
\ee
and the function
\bea
n_H(t,s)& =& \frac{\sqrt{\pi}}{2^{H-\frac{1}{2}}}I_{T,2,\frac{3-2H}{4}}^{H-\frac{1}{2}}(u^{H-\frac{1}{2}})I_{[0,t)}(s)\\\nonumber
&=& \frac{2^{1-H}\sqrt{\pi}}{\Gamma(H-\frac{1}{2})}s^{\frac{3}{2}-H}\int_0^t(x^2-s^2)^{H-\frac{3}{2}}dx \;I_{(0,t)}(s).\\\nonumber
\eea
The following theorem is due to Dzhaparidze and Van Zanten [5] (cf. Tudor [34]).
\vsp
\NI{\bf Theorem 2.2:} {\it The following representation holds, in distribution, for a sub-fBm $\zeta^H$:
\be
\zeta^H_t \stackrel {\Delta} = c_H\int_0^t n_H(t,s)dW_s, 0\leq t \leq T
\ee
where
\be
c_H^2= \frac{\Gamma (2H+1)\;\sin(\pi H)}{\pi}
\ee
and $\{W_t,t\geq 0\}$ is the standard Brownian motion.}
\vsp
Tudor [34] has defined integration of a non-random function $f(t)$ with respect to a sub-fBm $\zeta^H$ on an interval $[0,T]$ and obtained a representation of this integral as a Wiener integral for a suitable transformed function $\phi_f(t)$ depending on $H$ and $T.$ For details, see Theorem 3.2 in Tudor [34].
\vsp
Tudor [32] (cf. Tudor [34], p. 467) obtained the prediction formula for a sub-fBm. For any $0<H<1,$ and $0<a<t,$
\be
E[\zeta_t^H|\zeta_s^H, 0 \leq s \leq a]=\zeta_a^H+ \int_0^a \psi_{a,t}(u)d \zeta_u^H
\ee
where
\be
\psi_{a,t}(u)= \frac{2\; \sin(\pi (H-\frac{1}{2}))}{\pi} u(a^2-u^2)^{\frac{1}{2}-H}\int_a^t\frac{(z^2-a^2)^{H-\frac{1}{2}}}{z^2-u^2}z^{H-\frac{1}{2}}dz.
\ee
Let
\be
M_t^H=d_H \int_0^ts^{\frac{1}{2}-H}dW_s = \int_0^tk_H(t,s)d\zeta_s^H
\ee
where
\be
d_H=\frac{2^{H-\frac{1}{2}}}{c_H \Gamma(\frac{3}{2}-H)\sqrt{\pi}},
\ee
\be
k_H(t,s)=d_Hs^{\frac{1}{2}-H}\psi_H(t,s),
\ee
and
\ben
\psi_H(t,s)&=& \frac{s^{H-\frac{1}{2}}}{\Gamma(\frac{3}{2}-H)}[t^{H-\frac{3}{2}}(t^2-s^2)^{\frac{1}{2}-H}-\\\nonumber
&&\;\;\;\; (H-\frac{3}{2})\int_s^t(x^2-s^2)^{\frac{1}{2}-H}x^{H-\frac{3}{2}}dx]I_{(0,t)}(s).\\\nonumber
\een
\vsp
It can be shown that the process $ M^H=\{M^H_t, 0\leq t \leq T\}$ is a Gaussian martingale (cf. Tudor [34], Diedhiou et al. [4]) and is called the {\it sub-fractional fundamental martingale}. The filtration generated by this martingale is the same as the filtration $\{\clf_t, t\geq 0\}$ generated by the sub-fBm $\zeta^H$ and the quadratic variation $<M^H>_s$ of the martingale $M^H$ over the interval $[0,s]$ is equal to $w_s^H=\frac{d_H^2}{2-2H}s^{2-2H}= \lambda_H s^{2-2H}$ (say). For any measurable function $f:[0,T] \rightarrow R$ with $\int_0^Tf^2(s)s^{1-2H} ds<\ity,$ define the probability measure $Q_f$ by
\ben
\frac{dQ_f}{dP}|_{\clf_t}&=& \exp(\int_0^tf(s)dM_s^H-\frac{1}{2}\int_0^tf^2(s)d<M^H>(s))\\
&=&\exp(\int_0^tf(s)dM_s^H-\frac{d_H^2}{2}\int_0^tf^2(s)s^{1-2H}ds)
\een
where $P$ is the underlying probability measure. Let
\be
(\psi_H f)(s)=\frac{1}{\Gamma(\frac{3}{2}-H)}I_{0,2,\frac{1}{2}-H}^{H-\frac{1}{2}}f(s)
\ee
where, for $\alpha >0,$
\be
(I_{0,\sigma,\eta}^\alpha f)(s)=\frac{\sigma s^{-\sigma(\alpha+ \eta)}}{\Gamma(\alpha)}\int_0^s \frac{t^{\sigma(1+\eta)-1}f(t)}{(t^\sigma-s^\sigma)^{1-\alpha}}dt, s\in[0,T].
\ee
Then the following Girsanov type theorem holds for the sub-fBm process (Tudor [34]).
\vsp
\NI{\bf Theorem 2.3:} {\it The process
$$\zeta_t^H- \int_0^t(\psi_Hf)(s)ds, 0\leq t \leq T$$
is a sub-fbm with respect to the probability measure $Q_f.$ In particular, choosing the function $f\equiv a \in R$, it follows that the process $\{\zeta_t^H-at, 0\leq t \leq T\}$ is a sub-fBm under the probability measure $Q_f$ with $f \equiv a\in R.$}
\vsp
Let $Y=\{Y_t, t \geq 0\}$ be a stochastic process defined on the filtered probability space $(\Omega, \clf, (\clf_t, t\geq 0), P)$ and suppose the process $Y$ satisfies the stochastic differential equation
\be
dY_t= C(t) dt+ d\zeta_t^H, t \geq 0
\ee
where the process $\{C(t), t \geq 0\},$ adapted to the filtration $\{\clf_t, t \geq 0\},$ such that the process
\be
R_H(t)= \frac{d}{dw_t^H}\int_0^tk_H(t,s)C(s)ds, t \geq 0
\ee
is well-defined and the derivative is understood in the sense of absolute continuity with respect to the measure generated by the function $w_H.$ Differentiation with respect to $w_t^H$ is understood in the sense:
$${dw_t^H}= \lambda_H (2-2H)t^{1-2H}dt$$
and
$$\frac{df(t)}{dw_t^H}= \frac{df(t)}{dt} / \frac{dw_t^H}{dt}.$$
Suppose the process $\{R_H(t), 0\leq t \leq T\},$ defined over the interval $[0,T]$ belongs to the space $L^2([0,T],dw_t^H).$ Define
\be
\Lambda_H(t)=\exp\{\int_0^tR_H(s)dM_s^H-\frac{1}{2}\int_0^t[R_H(s)]^2dw_s^H\}
\ee
with $E[\Lambda_H(T)]=1$ and the distribution of the process $\{Y_t, 0\leq t \leq T\}$ with respect to the measure $P^Y= \Lambda_H(t)\;P$ coincides with the distribution of the process $\{\zeta_t^H, 0 \leq t \leq T\|$ with respect to the measure $P.$
\vsp
We call the process $\Lambda^H$ as the {\it likelihood process} or the Radon-Nikodym derivative $\frac{dP^Y}{dP}$ of the measure $P^Y$ with respect to the measure $P.$
\vsp
Tudor [34] derived the following Girsanov type formula.
\vsp
\NI{\bf Theorem 2.4:} {\it Suppose the assumptions of Theorem 2.2 hold. Define
\be
\Lambda_H(T)= \exp \{ \int_0^T R_H(t)dM_t^H - \frac{1}{2}\int_0^TR_H^2(t)d w_t^H\}.
\ee
Suppose that $E(\Lambda_H(T))=1.$ Then the measure $P^*= \Lambda_H(T) P$ is
a probability measure and the probability measure of the process $Y$ under
$P^*$ is the same as that of the process $V$ defined by}
\be
V_t= \int_0^t d\zeta_s^H, 0 \leq t \leq T.
\ee.
\newsection{Sub-Fractional  Vasicek model}
The following model was introduced by Vasicek [35] for modeling interest rates in finance. It is a model of the form

$$ dX_t=(\alpha-\beta X_t)dt+ \gamma dW_t, 0\leq t \leq T$$
where $\alpha,\beta,\gamma$ are positive real numbers and $\{W_t, t\geq 0\}$ is the standard Brownian motion. The parameter $\beta$ corresponds to the speed of recovery, the ratio $\alpha/\beta$ is the long-term average interest rate and the parameter $\gamma$ represents the stochastic volatility. The Vasicek model is used in finance, economics, biology, physics, chemistry, medicine , environmental studies and in other areas for modeling purposes. In a series of papers, Lohvinenko and Ralchenko [12,13,14,15] studied asymptotic properties of the maximum likelihood estimators of the parameters $\alpha$ and $\beta$ in the fractional Vasicek model

$$ dX_t=(\alpha-\beta X_t)dt+ \gamma dW_t^H, 0\leq t \leq T$$
as $T \raro \ity$ where $\alpha,\beta,\gamma$ are positive real numbers and $\{W_t^H, t\geq 0\}$ is the standard  fractional  Brownian motion with Hurst index $H >\frac{1}{2}$.

Our aim in this paper is to obtain the asymptotic properties of the maximum likelihood estimator for the parameters $\alpha, \beta$ in a  sub-fractional Vasicek model driven by a sub-fractional Brownian motion. For related results on estimation of parameters involved in processes driven by sub-fractional Brownian motion (mFBm), see Mendy [16], Xiao et  al. [37], Kuang and Liu [9], Kuang and Xie [10], Yu [41], Prakasa Rao [22,23,24,25,26,27,28] among others. 
\vsp
Let us consider the {\it sub-fractional Vasicek model} 
\be 
dX_t=(\alpha-\beta X_t)dt + d\zeta_t^H, t\geq 0 
\ee 
with $X_0=x_0 $ where $\alpha, \beta$ are unknown positive parameters with {\it known} Hurst index $H\in (1/2,1).$ In other words $X=\{X_t, t \geq 0\}$ is a stochastic process satisfying the stochastic integral equation 
\be X_t= x_0+\int_0^t (\alpha-\beta X_s)ds + \zeta_t^H, t \geq 0. 
\ee
The process $X=\{X_t, t \geq 0\}$ is termed as the {\it sub-fractional Vasicek process driven by a sub-fractional Brownian motion} also called the {\it sub-fractional Vasicek process}. This equation has a unique solution given by
\be
X_t=x_0e^{-\beta t}+\frac{\alpha}{\beta}(1-e^{-\beta t})+\int_0^te^{-\beta(t-s)}d\zeta_s^H, t\geq 0
\ee
where the integral 
$$\int_0^t e^{-\beta(t-s)}d\zeta_t^H$$
is interpreted as a Wiener integral with respect to a sub-fractional Brownian motion (cf. Tudor [34]).
Let
\be Q_H(t) =\frac{d}{d<M^H>_t}\int_0^tk_H(s,t) (\alpha-\beta X(s)) ds, t \geq 0. 
\ee
The process $\{Q_H(t), t \geq  0\}$ is well defined from the results in Tudor [34] and the sample paths of the process $\{Q_H(t), 0 \leq t \leq T \}$ belong almost surely to $L^2([0,T], d<M^H>_t).$ Define
\be 
Z_t= \int_0^t k_H(s,t)dX_s, t \geq 0. 
\ee 
Then the process
$Z= \{Z_t,t \geq 0\}$ is an $(\clf_t)$-semimartingale with the decomposition 
\be 
Z_t= \int_0^t Q_H(s)d<M^H>_s + M^H_t , t \geq 0
\ee
where $M^H$ is the fundamental Gaussian martingale. Let $P_\theta^T$ be the probability measure induced by the process $\{X_t, 0 \leq t \leq T\}$ when $\theta = (\alpha,\beta)$ is the true parameter. We write $Q_{H,\theta}(t) $ for $Q_H(t)$ hereafter  when $\theta$ is the true parameter. Following Theorem 2.4, we get that the Radon-Nikodym derivative of $P_\theta^T$ with respect to $P_0^T$ is given by
\be
L_T(\theta) \equiv \frac{dP_{\theta}^T}{dP_0^T}= \exp[ \int_0^T Q_{H,\theta}(s)dZ_s - \frac{1}{2}\int_0^T [Q_{H,\theta}(s)]^2d<M^H>_s].
\ee
\vsp
We now consider the problem of estimation of the parameter $\theta= (\alpha,\beta)$ based
on the observation of the process $ X= \{X_t , 0 \leq t \leq T\}$ or equivalently $\{Z_t, 0\leq t \leq T\}$ and
study the asymptotic properties of such estimators as $T \raro \ity.$ Let $\Theta= R_+^2.$
\vsp
\NI{Maximum Likelihood Estimation:}

Maximum likelihood estimator (MLE)  $\hat \theta_T$ is defined by the relation
\be
L_T(\hat \theta_T)= \sup_{\theta \in \Theta}L_T(\theta).
\ee
We assume that there exists a measurable MLE. Sufficient conditions can be given for the existence of such an estimator (cf. Lemma 3.1.2, Prakasa Rao [19]). Observe that
\bea
\Lambda_T(\theta) & = & \log L_T(\theta)\\\nnb 
&=& \int_0^T Q_{H,\theta}(s)dZ_s - \frac{1}{2} \int_0^T [Q_{H,\theta}(s)]^2d<M^H>_s\\\nnb
&=& \alpha Z_T-\beta \int_0^T P_H(t)dZ_t-\frac{1}{2}\alpha^2 <M^H>_T\\\nnb
&&\;\;\;\;+\alpha \beta\int_0^TP_H(t)d<M^H>_t-\frac{1}{2}\beta^2\int_0^T[P_H(t)]^2d<M^H>_t\\\nnb
\eea
where
\be
P_H(t)= \frac{d}{d<M^H>_t}\int_0^t k_H(s,t)X_s\;ds.
\ee
\vsp 
\NI{\bf Theorem 3.1:} {\it Suppose the parameter $\beta$ is known. Then the  MLE $\hat \alpha_T$ for $\alpha$ is
\be
\hat \alpha_T= \frac{Z_T+\beta\int_0^TP_H(t)d<M^H>_t}{<M^H>_T}. 
\ee
The MLE $\hat \alpha_T$  is unbiased, strongly consistent as $T \raro \ity.$ Furthermore the random variable $T^{1-H} (\hat \alpha_T-\alpha)$has a normal distribution with mean zero and variance $\lambda_H^{-1}$ depending on the Hurst index $H.$} 
\vsp
\NI{\bf Proof :} Maximizing the log-likelihood $\Lambda_T(\theta)$, lead to the equations
$$\frac{\partial \Lambda_T(\theta)}{\partial \alpha}= Z_T- \alpha <M^H>_T+\beta \int_0^TP_H(t)d<M^H>_t$$
and
$$\frac{\partial^2 \Lambda_T(\theta)}{\partial \alpha^2}= -<M^H>_T.$$
The equations given above imply that the MLE $\alpha_T$ of $\alpha$ is given by the equation (3.11). By Theorem 2.1, the process $Z$ has the representation
$$Z_T= \alpha<M^H>_T-\beta\int_0^TP_H(t)d<M^H>_t+M^H_T.$$
Applying this representation, it is easy to check that
$$\hat \alpha_T= \alpha+ \frac{M_T^H}{<M^H>_T}.$$
\vsp

Observe that the process $M^H$ is a martingale with the quadratic variation $<M^H>_t= \frac{d_H^2}{2-2H} t^{2-2H}=\lambda_H t^{2-2H} (say).$ 
Since $H <1,$ it follows that the function $<M^H>_T$ tends to infinity as $T \raro \ity.$ Hence , by the strong law of large numbers for martingales (cf. Liptser and Shiryayev [11], Theorem 2.6.10, Prakasa Rao [20]), it follows that
$$\frac{M_T^H}{<M^H>_T} \raro  0$$
almost surely as $T \raro \ity.$ Hence $\hat \alpha_T \raro \alpha$ almost surely as $T \raro \ity.$ Since the process $M^H$ is a Gaussian martingale with the quadratic variance $<M^H>,$ it follows that the random variable
$$\frac{M^H_T}{\sqrt{<M^H>_T}}$$
has the standard normal distribution for any fixed $T>0.$ . This in turn proves that the random variable
$$\lambda_H^{1/2}T^{1-H}(\hat \alpha_T-\alpha)$$
has the standard normal distribution. Hence
$$T^{1-H}(\hat \alpha_T-\alpha)$$
has the normal distribution with mean zero and variance $\lambda_H^{-1}.$
\vsp
\NI{\bf Theorem 3.2:} {\it Suppose the parameter $\alpha$ is known. Suppose that
$$\int_0^T[P_H(t)]^2d<M^H>_t \;\raro\; \ity$$
in probability as $T\raro \ity$ where $P_H(t)$ is as defined by the equation (3.10). Then the  MLE $\hat \beta_T$ for the parameter $\beta$ is
\be
\hat \beta_T = \frac{\alpha \int_0^T P_H(t)d<M^H>_t - \int_0^T P_H(t)dZ_t}{\int_0^T [P_H(t)]^2d<M^H>_t}. 
\ee
The estimator $\hat \beta_T$ is strongly consistent as $T \raro \ity.$ Furthermore the random variable
$$\sqrt{\int_0^T[P_H(t)]^2d<M^H>_t}(\hat \beta_T-\beta)$$
is asymptotically standard normal as $T\raro \ity.$}
\vsp
\NI{\bf Proof:} Maximizing the log-likelihood $\Lambda_T(\theta)$, lead to the equations
$$\frac{\partial \Lambda_T(\theta)}{\partial \beta}= -\int_0^TP_H(t)dZ_T+ \alpha \int_0^T P_H(t)d<M^H>_t-\beta \int_0^T[P_H(t)]^2d<M^H>_t$$
and
$$\frac{\partial^2 \Lambda_T(\theta)}{\partial \beta^2}= -\int_0^T[P_H(t)]^2d<M^H>_T.$$
which proves that the MLE $\hat \beta$ is given by the equation (3.12).
Note that
\be
dZ_t= \alpha d<M^H>_t -\beta P_H(t)d<M^H>_t+dM_t^H
\ee
and
\be
\int_0^TP_H(t)dZ_t= \alpha \int_0^T P_H(t)d<M^H>_t-\beta \int_0^T[P_H(t)]^2d<M^H>_t + \int_0^TP_H(t)d<M^H>_t.
\ee
Hence
\be
\hat \beta_T -\beta = \frac{\int_0^TP_H(t)dM_t^H}{\int_0^T[P_H(t)]^2d<M^H>_t}.
\ee
Since the process $M^H$ is a martingale with quadratic variation $<M^H>,$ the process 
$$\{\int_0^TP_H(t)dM^H_t, T \geq 0\}$$
is a local martingale with the quadratic variation 
$$\{\int_0^T[P_H(t)]^2d<M^H>_t, T \geq 0\}.$$
Observe that the process 
$$\{\int_0^T[P_H(t)]^2d<M^H>_t, T \geq 0\}.$$
is monotone increasing to infinity in probability as $T \raro \ity.$ Applying the strong law of large numbers for local martingales (cf. Liptser and Shiryayev [11], Theorem 2.6.10, Prakasa Rao [20]), it follows that $\hat \beta_T$ converges almost surely to $\beta$ as $T \raro \ity.$ Furthermore
\be
\sqrt{\int_0^T[P_H(t)]^2d<M^H>_t}(\hat \beta_T-\beta)= -\frac{\int_0^T P_H(t)dM^H_t}{\sqrt{\int_0^T[P_H(t)]^2d<M^H>_t}}
\ee
and the term on the right side of the above equation tends to the standard normal distribution as $T\raro \ity$ by the central limit theorem for local martingales (cf. Prakasa Rao [20]).
\vsp
\NI{\bf Theorem 3.3:} {\it Suppose both the parameters $\alpha$ and $\beta$ are unknown.  Then the MLEs of $\alpha$ and $\beta$ are given by
\be
\tilde \alpha_T= \frac{\int_0^TP_H(t)dZ_t\int_0^TP_H(t)d<M^H>_t-Z_T\int_0^T[P_H(t)]^2d<M^H>_t}{[\int_0^TP_H(t)d<M^H>_t]^2-<M^H>_T\int_0^T[P_H(t)]^2d<M^H>_t} 
\ee
and}
\be
\tilde \beta_T= \frac{<M^H>_T\int_0^TP_H(t)dZ_t -Z_T\int_0^TP_H(t)d<M^H>_t}{[\int_0^TP_H(t)d<M^H>_t]^2-<M^H>_T\int_0^T[P_H(t)]^2d<M^H>_t}. 
\ee
\vsp
\NI{\bf Proof:} Maximizing the log-likelihood $\Lambda_T(\theta)$ with respect to the parameter $\alpha$ and $\beta$ simultaneously lead to the equations
\be
\frac{\partial \Lambda_T(\theta)}{\partial \alpha}= Z_T-\alpha<M^H>_T+\beta\int_0^TP_H(t)d<M^H>_t=0
\ee
and
\bea
\frac{\partial \Lambda_T(\theta)}{\partial \beta}\\\nnb
= -\int_0^TP_H(t)dZ_T+\alpha \int_0^TP_H(t)d<M^H>_t-\beta\int_0^T[P_H(t)]^2d<M^H>_t=0.
\eea
Solving these equations, we obtain the estimators $\tilde \alpha_T$ and $\tilde \beta_T$ as given by the equations (3.17) and (3.18) respectively. Observe that
$$\frac{\partial^2 \Lambda_T(\theta)}{\partial \alpha^2}=-\frac{<M^H>_T}{T}<0,$$
$$\frac{\partial^2 \Lambda_T(\theta)}{\partial \beta^2}= -\int_0^T[P_H(t)]^2d<M^H>_t<0,$$
and
$$ \frac{\partial^2 \Lambda_T(\theta)}{\partial \alpha^2} \frac{\partial^2 \Lambda_T(\theta)}{\partial \beta^2}-[\frac{\partial^2 \Lambda_T(\theta)}{\partial \alpha \partial \beta}]^2= <M^H>_T\int_0^T[P_H(t)]^2d<M^H>_t-[\int_0^TP_H(t)d<M^H>_t]^2<0$$
by the Cauchy-Schwartz inequality which implies that the estimators $\tilde \alpha_T$ and $\tilde \beta_T$ maximize the likelihood and hence are the MLEs of $\alpha $ and $\beta $ respectively. An application of the representation of the process $\{Z_t, 0\leq t \leq T$ given by Theorem 2.1 implies that
\be
\tilde \alpha_T-\alpha= \frac{\int_0^TP_H(t)dM^H_t\int_0^TP_H(t)d<M^H>_t-M^H_T\int_0^T[P_H(t)]^2d<M^H>_t}{[\int_0^TP_H(t)d<M^H>_t]^2-M^H_T\int_0^T(P_H(t))^2d<M^H>_t}
\ee
and
\be
\tilde \beta_T-\beta = \frac{<M^H>_T \int_0^T P_H(t)dM^H_t-M^H_T\int_0^TP_H(t) d<M^H>_t}{[\int_0^TP_H(t)d<M^H>_t]^2-M^H_T\int_0^T(P_H(t))^2d<M^H>_t}.
\ee
\vsp
\NI{\bf Remarks:} Even though the form of the function $k_H(t,s)$ is known, due to its complicated form, it is not possible to use the methods in Lohvinenko and Ralchenko [12,13,14,15] to study the asymptotic distribution of $(\tilde \alpha_T,\tilde \beta_T)$ or the asymptotic marginal distributions of $\tilde \alpha_T$ and $\tilde \beta_T$ after suitable scaling as $T \raro \ity.$ However, following ideas in Lohvinenko and Ralchenko [13], we will transform the problem to the study of maximum likelihood estimation for the parameters of sub-fractional Vasicek model to that of estimation of parameters for a sub-fractional Ornstein-Uhlenbeck process (Mendy [16], Yu [41], Es-Sebaiy-Es-sebaiy [7], Xiao et al. [37] ) and derive the asymptotic properties of the corresponding MLE. We consider the case $H>\frac{1}{2}.$
\vsp
\newsection{Alternate approach}
Consider the following process
\be
U_t= \int_0^te^{-\beta(t-s)}d\zeta_t^H, t \geq 0.
\ee
Then the process $\{U_t, t \geq 0\}$ is a sub-fractional Ornstein-Uhlenbeck process  and it is the solution of the equation
\be
dU_t= -\beta U_t dt+ d\zeta_t^H, U_0=0
\ee
The sub-fractional Vasicek model defined by the stochastic differential equation
\be
dX_t= (\alpha-\beta X_t)dt+ d\zeta_t^H, t \geq 0, X_0=x_0
\ee
can be rewritten in the form
\be
X_t= \frac{\alpha}{\beta}+ (x_0-\frac{\alpha}{\beta})e^{-\beta t}+ U_t. 
\ee
Observe that
\bea
\;\;\;\;\\\nnb
P_H(t)&=& \frac{d}{d<M^H>_t}\int_0^tk_H(t,s)X_sds\\\nnb
&=& \frac{d}{d<M^H>_t}\int_0^tk_H(t,s)[\frac{\alpha}{\beta}+ (x_0-\frac{\alpha}{\beta})e^{-\beta s}+U_s]ds\\\nnb
&=& \frac{\alpha}{\beta}\frac{d}{d<M^H>_t}\int_0^tk_H(t,s)ds + (x_0-\frac{\alpha}{\beta})\frac{d}{d<M^H>_t}\int_0^tk_H(t,s) e^{-\beta s}ds\\\nnb
&&\;\;\;\;+\frac{d}{d<M^H>_t}\int_0^tk_H(t,s)U_sds\\\nnb
&=& \frac{\alpha}{\beta}J(t)+(x_0- \frac{\alpha}{\beta})V_H(t)+\tilde P_H(t)
\eea
where
\be
J(t)= \frac{d}{d<M^H>_t}\int_0^tk_H(t,s)ds,
\ee
\be
\tilde P_H(t)=\frac{d}{d<M^H>_t}\int_0^tk_H(t,s)U_sds
\ee
and
\be
V_H(t)=  \frac{d}{d<M^H>_t}\int_0^tk_H(t,s) e^{-\beta s}ds.
\ee
Suppose that
$$ \frac{1}{T}\int_0^T[\tilde P_H(t)]^2d<M^H>_t \raro C_{H,\beta}$$
in probability as $T\raro \ity$ for some positive constant $C_{H,\beta}$. This in turn implies that 
$$ \int_0^T[\tilde P_H(t)]^2d<M^H>_t \raro \ity$$
in probability as $T \raro \ity$ and 
$$ \frac{1}{\sqrt{T}}\int_0^T\tilde P_H(t)dM^H_t \raro N(0, C_{H,\beta})$$
in distribution as $T \raro \ity$ by the central limit theorem for local martingales (cf. Prakasa Rao (1999b)) for some positive constant $C_{H,\beta}$ depending on $H$ and $\beta.$. Let $\beta_T^*$ be the maximum likelihood estimator of $\beta.$ It can be checked that
\be
\beta_T^*-\beta= \frac{\int_0^T\tilde P_H(t)dM^H_t}{\int_0^T[\tilde P_H(t)]^2d<M^H>_t}.
\ee
Furthermore
\be
\sqrt{T}(\beta_T^*-\beta) \raro N(0,C_{H,\beta}^{-1})
\ee
in distribution as $T\raro \ity.$ 
\vsp
\NI{\bf Theorem 4.1:} Suppose that  
$$ \frac{1}{T}\int_0^T[\tilde P_H(t)]^2d<M^H>_t \raro C_{H,\beta}$$
in probability as $T\raro \ity$ for some positive constant $C_{H,\beta}$ where $\tilde P_H(t)$ is as defined by the equation (3.29). Let $\beta_T^*$ be the maximum likelihood estimator of $\beta.$ Then
\be
\sqrt{T}(\beta_T^*-\beta) \raro N(0,C_{H,\beta}^{-1})
\ee
in distribution as $T\raro \ity.$ Further suppose that 
$$\lim_{T\raro \ity}E|\frac{1}{T}\int_0^T[\tilde P_H(t)]^2d<M^H>_t|^{p}<\ity, p\geq 1.$$
Then
\be 
E[(\sqrt{T}(\beta^*_T-\beta))^p] \raro E[(\sqrt{C_{H,\beta}^{-1}}Z)^p]
\ee
as $T \raro \ity$ where $Z$ is N(0,1) holds for all integers $p \geq 1.$
\vsp
\NI{\bf Proof:} Asymptotic normality of the estimator $\beta_T^*$ as $T \raro \ity $ was proved by the arguments given above. Note that
\bea
\;\;\;\;\\\nnb
\frac{1}{<M^H>_T}\int_0^TP_H(t)d<M^H>_t &=& \frac{1}{<M^H>_T}\int_0^T[\frac{\alpha}{\beta}J(T)
+(x_0-\frac{\alpha}{\beta})V_H(t)+\tilde P_H(t)]d<M^H>_t\\\nnb
&=& \frac{\alpha}{\beta} J(T)+ (x_0-\frac{\alpha}{\beta})\frac{1}{<M^H>_T} \int_0^TV_H(t) d<M^H>_t\\\nnb
&&\;\;\;\;+ \frac{1}{<M^H>_T} \int_0^T \tilde P_H(t) d<M^H>_t.\\\nnb
\eea
Note that the convergence of moments 
\be 
E[(\sqrt{T}(\beta^*_T-\beta))^p] \raro E[(\sqrt{C_{H,\beta}^{-1}}Z)^p]
\ee
as $T \raro \ity$ where $Z$ is N(0,1) holds for all integers $p \geq 1$ if the family $(\sqrt{T}(\hat\beta_T-\beta))^p$ is uniformly integrable over $T$ for all integers $p\geq 1.$ Observe that
\bea
\;\;\;\;\\\nnb
 (E[|\sqrt{T}(\beta^*_T-\beta)|^p])^2 &\leq & E|\frac{1}{T}\int_0^T[\tilde P_H(t)]^2d<M^H>_t|^{-2p}E|\frac{1}{\sqrt{T}}\int_0^T\tilde P_H(t)dM^H_t|^{2p}\\\nnb
&\leq & E|\frac{1}{T}\int_0^T[\tilde P_H(t)]^2d<M^H>_t|^{-2p} C_pE|\frac{1}{T}\int_0^T[\tilde P_H(t)]^2d<M^H>_t|^{p}\\\nnb
\eea
where the last bound holds by the Burkholder-Davis-Gundy inequality with an absolute constant $C_p.$ Hence the limit in the equation (4.14) holds by the de la Vallee Poussin theorem since
$$\lim_{T\raro \ity}E|\frac{1}{T}\int_0^T[\tilde P_H(t)]^2d<M^H>_t|^{p}<\ity$$
for all integers $p \geq 1$ by hypothesis.
\vsp 
\NI{\bf Acknowledgment:} Work in this paper was supported under the scheme `INSA Senior Scientist" at the CR Rao Advanced Institute for Mathematics, Statistics and Computer Science, Hyderabad, India.
\vsp
\NI {\bf References :}
\begin{description}

\item {[1]} T. Bojdecki, T., L.G. Gorostiza and  A. Talarczyk, Sub-fractional Brownian motion and its relation to occupation times, {\it Statist. Probab. Lett.}, {\bf 69} (2004) no.4, 405-419.

\item {[2]} A. Chronopoulu, A. and F,G. Viens, Estimation and pricing under long-memory stochastic volatility, {\it Ann. finance}, {\bf 8}(2012) , 379-403.

\item {[3]} S. Corlay, S., J. Lobovits, J., and J.L. Vehel, Multifractional stochastic volatility models, {\it Math. Finance}, {\bf 24}(2014) , 364-402.

\item {[4]} A. Diedhiou, A., C. Manga, C. and I. Mendy, Parametric estimation for SDEs with additive sub-fractional Brownian motion, {\it J. Numer. Math. and Stoch.}, {\bf 3} (2011), no.1, 37-45.

\item {[5]} K. Dzhaparidze, and H. Van Zanten, A series expansion of fractional Brownian motion, {\it Probab. Theory Related Fields}, {\bf 130} (2004), no.1, 39-55.

\item {[6]} M. El Machkouri, K. Es-Sebaiy and Y. Ouknine, Least squares estimation for non-ergodic Ornstein-Uhlenbeck processes driven by Gaussian processes, {\it J. Korean Statist. Soc.} {\bf 45} (2016), no.3, 329-341.

\item {[7]} K. Es-Sebaiy, K. and M. Es-Sebaiy, Estimating drift parameters in a non-ergodic fractional Vasicek model, arXiv 1909.06155v3 [math Pr] 9 May 2020.

\item {[8]} R. Hao, Y. Liu, and S. Wang, Pricing credit default swap under fractional Vasicek interest rate model, {\it J. Math. Finance}, {\bf 4} (2014) , 10-20.  

\item {[9]} N. Kuang and B. Liu, Parameter estimations for the sub-fractional Brownian motion with drift at discrete observation, {\it Brazilian Journal of Probability and Statistics}, {\bf 29} (2015), no.4, 778-789.

\item {[10]} N. Kuang and H. Xie, Maximum likelihood estimator for the sub-fractional Brownian motion approximated by a random walk, {\it Ann. Inst. Statist. Math.}, {\bf 67} (2015), no.1, 75-91.

\item {[11]} R.S. Liptser and A.N. Shiryayev, {\it The Theory of Martingales}, Kluwer, Dordrecht, 1989.

\item {[12]} S. Lohvinenko and K. Ralchenko, Asymptotic properties of parameter estimators in  fractional Vasicek model, {\it Lithuanian J. Statist.}, {\bf 55} (2016), 102-111.

\item {[13]} S. Lohvinenko and K. Ralchenko, Maximum likelihood estimation in the fractional Vasicek model, {\it Lithuanian J. Statist.}, {\bf 56} (2017) , 77-87.

\item {[14]} S. Lohvinenko, and K. Ralchenko, Asymptotic distribution of the maximum likelihood estimator in the fractional Vasicek model, {\it Theor. Probability and Math. Statist.}, {\bf 99} (2018), 134-151.

\item {[15]} S. Lohvinenko and K. Ralchenko, Maximum likelihood estimation in the non-ergodic fractional Vasicek model, {\it Modern Stochastics: Theory and Applications}, {\bf 6} (2019), 377-395.

\item {[16]} I. Mendy, Parametric estimation for sub-fractional Ornstein-Uhlenbeck process, {\it J. Stat. Plan. Infer.}, {\bf 143} (2013), 663-674.

\item {[17]} Y. Mishura and M. Zili, {\it Stochastic Analysis of Mixed Fractional Gaussian Processes}, ISTE Press, London, 2018.

\item {[18]} Y. Mishura, {\it Stochastic Calculus for Fractional Brownian Motion and Related Processes}, Lecture Notes in Math. 1929, Springer, Berlin, 2008.

\item {[19]} B.L.S. Prakasa Rao, {\it Asymptotic Theory of Statistical Inference}, Wiley, New York, 1987.

\item {[20]} B.L.S. Prakasa Rao, {\it Semimartingales and Their Statistical Inference}, CRC Press, Boca Raton and Chapman and Hall, London, 1999.

\item {[21]} B.L.S. Prakasa Rao, {\it Statistical Inference for Fractional Diffusion Processes}, Wiley, London, 2010.

\item {[22]} B.L.S. Prakasa Rao,  On some maximal and integral inequalities for sub-fractional Brownian motion, {\it Stochastic Anal. Appl.}, {\bf 35} (2017), no.2,  279-287.

\item {[23]} B.L.S. Prakasa Rao, Optimal estimation of a signal perturbed by a sub-fractional Brownian motion, {\it Stochastic Anal. Appl.}, {\bf 35} (2017), no.3, 533-541.

\item {[24]} B.L.S. Prakasa Rao, Parameter estimation for linear stochastic differential equations driven by sub-fractional Brownian motion, {\it Random Oper. and Stoch. Equ.}, {\bf 25} (2017), no.4, 235-247.

\item {[25]} B.L.S. Prakasa Rao, Berry-Esseen type bound for Fractional Ornstein-Uhlenbeck type process driven by sub-fractional Brownian motion, {\it Theory Stoch. Proc.}, {\bf 23} (2018), no.1, 82-92.

\item {[26]} B.L.S. Prakasa Rao, More on maximal inequalities for sub-fractional Brownian motion, {\it Stochastic Anal. Appl.}, {\bf 38} (2020) ,no.2,  238-247.

\item {[27]} B.L.S. Prakasa Rao, Nonparametric estimation of trend for stochastic differential equations driven by sub-fractional Brownian motion, {\it Random Oper. Stoch. Equ.}, {\bf 28} (2020) , 113-122.

\item {[28]} B.L.S. Prakasa Rao, Maximum likelihood estimation in the mixed fractional Vasicek model, Technical Report, CR Rao Advanced Institute of Mathematics, Statistics and Computer Science, Hyderabad, India (2020).

\item {[29]} G.J. Shen and L.T. Yan, Estimators for the drift of sub-fractional Brownian motion, {\it Communications in Statistics-Theory and Methods}, {\bf 43} (2014), no.8, 1601-1612.

\item {[30]} L. Song and K. Li, Pricing option with stochastic interest rates and transaction costs in fractional Brownian markets, {\it Discrete Dyn. Nat. Soc.}, (2018) 7056734-8.

\item {[31]} C. Tudor, Some properties of the sub-fractional Brownian motion, {\it Stochastics}, {\bf 79} (2007), no.5, 431-448.

\item {[32]} C. Tudor, Prediction and linear filtering with sub-fractional Brownian motion, Preprint (2007).

\item {[33]} C. Tudor, Some aspects of stochastic calculus  for the sub-fractional Brownian motion, {\it An. Univ.  Bucaresti, Mat.}, {\bf 57} (2008), no.2,  199-230.

\item {[34]} C. Tudor, On the Wiener integral with respect to a sub-fractional Brownian motion on an interval, {\it J. Math. Anal. Appl.}, {\bf 351} (2009),no.1, 456-468.

\item {[35]} O. Vasicek, An equilibrium characterization of the term structure. {\it J. finance Econ.}, {\bf 5} (1977), 177-188.

\item {[36]} W. Xiao and J. Yu, Asymptotic theory for rough fractional Vasicek models. {\it Econ. Lett.}, {\bf 177} (2019), 26-29.

\item {[37]} W. Xiao and J. Yu, Asymptotic theory for estimating drift parameters in the fractional Vasicek model, {\it Econom. Theory}, {\bf 35} (2019) , 198-231.

\item {[38]} W. Xiao, W. Zhang, X. Zhang and X. Chen, The valuation of equity warrants under the fractional Vasicek process of the short-term interest rate. {\it Physica A}, {\bf 394} (2014), 320-337.

\item {[39]} W. Xiao, X. Zhang, and Y. Zuo, Least squares estimation for the drift parameters in the sub-fractional Vasicek processes, {\it J. Stat. Plan. Inf.}, {\bf 197} (2018), 141-155.

\item {[40]} L. Yan, L., G. Shen and K. He, Ito's formula for a sub-fractional Brownian motion, {\it Commun. Stoch. Anal.}, {\bf 5} (2011), no.1, 135-159.

\item {[41]} Q. Yu, Statistical inference for Vasicek-type model driven by self-similar Gaussian processes, {\it Comm. in Statist.-Theory and Methods}, {\bf 49} (2020), 471-484.

\end{description}
\end{document}